\documentclass{article} 

\usepackage{epsfig,amsmath}
\usepackage[psamsfonts]{amssymb}
\usepackage[latin1]{inputenc}

\newcommand{\wt}{\widetilde}
\newcommand{\wh}{\widehat}

\newtheorem{proposition}{\bf Proposition}
\newtheorem{definition}{\bf Definition}
\newtheorem{theorem}{\bf Theorem}
\newtheorem{lemma}{\bf Lemma}
\newtheorem{corollary}{\bf Corollary}

\def\C{{\mathbb C}}
\def\D{{\mathbb D}}

\def\Q{{\mathbb Q}}
\def\Z{{\mathbb Z}}
\def\R{{\mathbb R}}

\def\N{{\mathbb N}}
\def\U{{\mathbb U}}

\def\a{\alpha}
\def\e{\varepsilon}
\def\cal{\mathcal}
\def\re{{\rm Re}}
\def\im{{\rm Im}}
\def\rad{{\rm rad}}

\def\proof{\noindent{\bf Proof. }}
\def\remark{\vskip.2cm \noindent{\bf Remark. }}
\def\qedprop{\hfill{\hbox{%
  \hskip 1pt%
  \vrule width 7pt height 6pt depth 1.5pt%
  \hskip 1pt}} \vskip.3cm}
\def\qedlem{\hfill{$\square$} \vskip.3cm}

\def\ds{\displaystyle}

\begin{document}

\title{On the Size of Quadratic Siegel Disks: Part I.}
\author{Xavier Buff\thanks{
Laboratoire Emile Picard, UMR CNRS 5580,
UFR MIG, Universit{\'e} Paul Sabatier,
118 route de Narbonne,
31062 Toulouse Cedex 4,
France.
\texttt{buff@picard.ups-tlse.fr; http://picard.ups-tlse.fr/\~{
}buff}. }
and Arnaud Chéritat\thanks{
Laboratoire Emile Picard, UMR CNRS 5580,
UFR MIG, Universit{\'e} Paul Sabatier,
118 route de Narbonne,
31062 Toulouse Cedex 4,
France.
\texttt{cheritat@picard.ups-tlse.fr; http://picard.ups-tlse.fr/\~{
}cheritat}. }}

\date{May 2003}

\maketitle

\abstract{
If $\a$ is an irrational number, we let $\{p_n/q_n\}_{n\geq 0}$, be the
approximants given by its continued fraction expansion. The Bruno series
$B(\a)$ is defined as
$$B(\a)=\sum_{n\geq 0} \frac{\log q_{n+1}}{q_n}.$$

The quadratic polynomial $P_\a:z\mapsto e^{2i\pi \a}z+z^2$ has an
indifferent fixed       
point at the origin. If $P_\a$ is linearizable, we let $r(\a)$ be the
conformal radius of the Siegel disk and we set $r(\a)=0$ otherwise.
Yoccoz proved that if $B(\a)=\infty$, then $r(\a)=0$ and $P_\a$ is
not linearizable. 
In this article, we present a different proof and we show that 
there exists a constant $C$ such that
for all irrational number $\a$ with $B(\a)<\infty$, we have 
$$B(\a)+\log r(\a) < C.$$ 
Together with former results of Yoccoz (see \cite{y}), this proves
the conjectured boundedness of $B(\a)+\log r(\a)$.
}

\section{Introduction.}

In this article, we are interested in the dynamics of quadratic
polynomials $P_\a:z\mapsto e^{2i\pi\a}z+z^2$, $\a\in \C$. 
When $\a$ is real, the quadratic polynomial $P_\a$ has an
indifferent fixed point at $0$ and it is \em linearizable \em if it is
conjugate to the rotation $z\mapsto e^{2i\pi\a}z$ in a
neighborhood of $0$. The arithmetic nature of $\a$ will play a
central role. We denote by $\{p_n/q_n\}_{n\geq 0}$ the approximants to
$\a$ given by its continued fraction expansion (see
appendix~\ref{app_cf}).

\remark Every time we use the notation $p/q$ for a rational
number, we mean that $q>0$ and $p$ and $q$ are coprime.
\vskip.2cm 

The first result of linearizability is due to C.L. Siegel \cite{si} in
1942.
He proved that when 
$$\log q_{n+1} = {\cal O}(\log q_n),$$
every germ $z\mapsto e^{2i\pi \a}z+{\cal O}(z^2)$ is
linearizable. Around 1965, following Siegel's ideas, Bruno \cite{bru}
proved that every germ $z\mapsto e^{2i\pi \a}z+{\cal O}(z^2)$ is
linearizable under the weaker assumption:
$$\sum_{n=0}^{+\infty}\frac{\log q_{n+1}}{q_n}<+\infty.$$
An irrational number $\a$ satisfying this condition is called a Bruno
number, and the sum on the left-hand side of the inequality is
noted $B(\a)$.
In 1987, Yoccoz \cite{y} has completely solved the problem,
showing that when $\a\in \R$ is not a Bruno number, the quadratic
polynomial $P_\a$ is not linearizable. Yoccoz first proved that if
there were a non linearizable germ $z\mapsto e^{2i\pi \a}z+{\cal O}(z^2)$,
then the quadratic polynomial $P_\a$ would not be  
linearizable. He then proved the existence of such a germ. 

More recently, in \cite{c}, the second author proved directly the non
linearizability of $P_\a$ when $\a$ is not a Bruno number. His
proof consists in proving that when $\a$ is not a Bruno number,
$0$ is accumulated by periodic points of $P_\a$ (P\'erez-Marco \cite{pm} in
fact proved that not only periodic points, but whole cycles, accumulate $0$). 

Let us now introduce the notion of conformal radius. 

\begin{definition}
If $U\subset \C$ is a hyperbolic domain containing $0$, the conformal radius
$\rad(U)$ of $U$ at $0$ is equal to 
$|\pi'(0)|$ where $\pi:(\D,0)\to (U,0)$ is a universal covering.
\end{definition}

\remark When $U$ is simply connected, for example in the case of a Siegel
disk, $\pi$ is the Riemann mapping, and $\rad(U)$
is the classical conformal radius. 
\vskip.2cm

Assume that $0\in V\subset U$. Then the universal covering
$\pi_V:(\D,0)\to (V,0)$ lifts to a mapping $f:(\D,0)\to (\D,0)$ such
that $\pi_V=\pi_U\circ f$. By Schwarz's lemma, we have $|f'(0)|\leq 1$,
and thus, $\rad(V) = |\pi'_V(0)|\leq |\pi'_U(0)|=\rad(U)$:
\[\rad(V) \leq \rad(U).\]

The work by Yoccoz \cite{y} already
provides a control of the conformal radius
$r(\a)$ of the Siegel disk of $P_\a$.
First, if $B(\a)<\infty$, then 
$$-B(\a)+C \leq \log r(\a)$$ 
for some universal constant $C$. Second, there exists a function
$C'(\e)$, such that
$\forall \e>0$,
$$\log r(\a) \leq -(1-\e)B(\a) + C'(\e).$$
In other words, $$ C \leq \log r(\a) + B(\a) \leq \e B(\a) +
  C'(\e).$$
In this article, we will prove the following results. 

\begin{theorem}
Assume $\a\in \R\setminus\Q$ is an irrational number and let $p_n/q_n$ be the
approximants to $\a$. For $N\geq 0$, 
let $V_N$ be the complement in $\C$ of the external ray
of $P_\a$ of argument $0$ and the periodic points of period less 
than or equal to $q_N$. Then, 
$$\log\rad(V_N) +\sum_{n=0}^N \frac{\log q_{n+1}}{q_n} < 16.$$
\end{theorem}

If $P_\a$ has a Siegel disk, it must be contained in the intersection of those
sets $V_N$. Thus, we have the following corollary. 

\begin{corollary}
If $B(\a)=\infty$, then $r(\a)=0$.
If $B(\a)<\infty$, then 
$$B(\a)+\log r(\a) < 16.$$
\end{corollary}

Note that we did not try to get the best possible
constants. The proof we give is a quantification of the proof given in \cite{c}
(with some minor modifications), together with a big improvement in
one inequality.

\section{Sketch of the proof.}

Let us first present the main steps of the proof. Note that the functions
$\a\mapsto B(\a)$ and $\a\mapsto r(\a)$ are even and periodic of
period $1$: for $B(\a)$ it is proved in appendix~\ref{app_cf}, and
for $r(\a)$, the periodicity comes from $P_{\a+1} = P_\a$, and 
the other claim from the fact that $P_{-\a}$ and $P_\a$ are
conjugated by an isometry (namely $z \mapsto \bar{z}$).
Thus, without loss of generality, we may assume
$\a\in ]0,1/2[\setminus 
\Q$.

\vskip.2cm
\noindent{\bf Step 1.}
For each rational number $p/q$, the polynomial $P_{p/q}$ has a
parabolic fixed point at $0$. When $\a$ is sufficiently close to
$p/q$, this parabolic point splits into a simple fixed point at $0$
and a periodic cycle of period $q$ which is close to $0$. The first
step consists in studying the dependence of this
cycle on $\a\in \C$. Roughly speaking, as long as the cycle does
not collide with another cycle, it is possible to follow it
holomorphically. More precisely, we have the following two
statements. The proofs are given in section \ref{sect3} below. 

\begin{definition}
For each rational number $p/q$, let $R(p/q)$ be the largest real
number such that $P_\a^{\circ q}$ has no multiple
fixed point for $\a \in B\big(\frac{p}{q},R(\frac{p}{q})\big)
\setminus \big\{\frac{p}{q}\big\}$.
Moreover, set $r(p/q)=[R(p/q)]^{1/q}$. 
\end{definition}

The proofs of the two following propositions are
detailed in \cite{c} and \cite{bc}, but for completeness, we sketch
them in section~\ref{sect3}.

\begin{proposition}\label{functionchi}
Let $p/q$ be a rational number, and $\zeta= e^{2i\pi p/q}$.
There exists an analytic function
$\chi:B(0,r(p/q))\to \C$ such that 
$\chi(0)=0$ and for any $\delta\in
B(0,r(p/q))\setminus\{0\}$, $\chi(\delta)\neq 0$ and the set 
$$\Big< \chi(\delta),\chi(\zeta\delta), 
\chi(\zeta^2\delta),\ldots,\chi(\zeta^{q-1}\delta)\Big>$$
forms a cycle of period $q$ of $P_{p/q+\delta^q}$.
We will note $\chi = \chi_{p/q}$, since it depends on $p/q$.
\end{proposition}

The proof is a simple application of the implicit function theorem.

\begin{proposition}[\sc Key inequality]\label{constantC}
For any rational number $p/q$, we have 
$$R(p/q)\geq \frac{1}{q^3}.$$
\end{proposition}
The proof relies on the Yoccoz inequality and on a
combinatorial theorem.
Note it is probably not optimal~:
the correct order is conjectured to be equal to $2$.

\vskip.2cm
\noindent{\bf Step 2.}
The polynomial $P_\a$ is a monic polynomial. It is affinely conjugate to
the polynomial $z\mapsto z^2+c$ with $c=e^{2i\pi\a}/2-e^{4i\pi\a}/4$.
As long as $c\notin [1/4,+\infty[$, there is a well-defined external ray
${\cal R}_0(\a)$ of argument $0$ which does not bifurcate and lands at a
repelling fixed point.
Note that $c\in[1/4,+\infty[$ $\Longleftrightarrow$ $\re(e^{2i\pi \a})=1$.

\begin{definition}
Denote by ${\cal B}$ the set of parameters $\a\in \C$ such that
$\re(e^{2i\pi \a})=1$. For each rational number $p/q$ with $q\geq 2$, let
$R'(p/q)$ be the largest real number such that 
$$B\left(\frac{p}{q},R'(p/q)\right)\subset \C\setminus {\cal B}.$$
\end{definition}

\begin{proposition}
For each rational number $p/q$ with $q\geq 2$, when $\a$ ranges in
$B(p/q,R'(p/q))$, the external ray ${\cal R}_0(\a)$ does
not bifurcate and lands at a repelling fixed point located
at $1-e^{2i\pi \a}$.  
\end{proposition}

\proof Since ${\cal R}_0(\a)$ does not bifurcate, it moves
holomorphically together with its landing point. This landing point must be a
fixed point of $P_\a$ and it cannot be $0$ since for $\a=p/q$,
and $q \geq 2$, the
external rays landing at $0$ are not fixed whereas ${\cal R}_0(\a)$ is
fixed. 
\qedprop

\begin{proposition}\label{constantC'}
For any rational number $p/q$ with $q\geq 2$, we have
$$R'(p/q)\geq \frac{1}{q^2}.$$ 
\end{proposition}

The proof is given in section \ref{sect4} below.
Note that in particular, $R'(p/q) \geq 1/q^3$.

\vskip.2cm
\noindent{\bf Step 3.}
Let us now assume that $\a_0\in ]0,1/2[\setminus \Q$ is an irrational
number and let $\{p_n/q_n\}_{n\geq 0}$ be the approximants given by 
its continued fraction expansion.
Then, for $n\geq 0$, $q_n$ is bounded 
from below by the $n$-th Fibonacci number $F_n$ ($F_0=1$, $F_1=2$, 
$F_{n+1}=F_n+F_{n-1}$).

The next definition will be better understood if the reader keeps in mind
that, according to classical properties of approximants, for all $n\in\N$,
\[
  \frac{1}{2 q_n q_{n+1}} < \left| \a_0-\frac{p_n}{q_n} \right| < \frac{1}{q_n q_{n+1}}.
\]

\begin{definition}[\sc Good approximants]\label{goodapproximants}
Let ${\cal N}$ be the set of integers $n\geq 1$ such that 
$q_{n+1}>2q_n^2$. Let $\{n_i\}_{i\geq 1}$, be the sequence
of those integers $n$ ordered increasingly. 

For $i\geq 1$, let 
\begin{itemize}
\item $B_i$ be the disk centered at
$p_{n_i}/q_{n_i}$ with radius $1/q_{n_i}^3$,
\item $B_i^*$ be the punctured disk
$B_i\setminus \{p_{n_i}/q_{n_i}\}$,
\item $D_i$ be the disk centered at
$p_{n_i}/q_{n_i}$ with radius $1/q_{n_i}^2$ and 
\item $U_i$ be the disk centered at $0$ with radius 
$(1/q_{n_i}^3)^{1/q_{n_i}}$.
\end{itemize}
\end{definition}

\remark Note that the set ${\cal N}$ may be finite, or even empty,
for example if $\a=(3-\sqrt {5})/2$. 
\remark The choice of the condition $q_{n+1}>2q_n^2$ and of the
  radius $1/q_{n_i}^3$ of $B_i$ are related
to the term $1/q^3$ in proposition~\ref{constantC}.
\vskip.2cm

The set ${\cal N}$ has been chosen so that the following
two propositions hold. 

\begin{proposition}\label{goodapprox1}
We have $B_1\subset D_1\subset \C\setminus {\cal B}$
 and for all $i\geq 1$, 
$\a_0\in B_i$ and
$$B_{i+1}\subset D_{i+1}\subset B_i^*.$$
\end{proposition}

\begin{proposition}\label{goodapprox2}
Moreover, for any $N\geq 1$
$$\sum_{n=1}^{N}
\frac{\log q_{n+1}}{q_{n}} < 
\sum_{\underset{\scriptstyle n_i\leq N}{i\geq 1}} 
\frac{\log q_{1+n_i}}{q_{n_i}}+
\sum_{n\in [1,N]\setminus {\cal N}} 
\frac{\log 2F_{n}^2}{F_{n}}.$$
\end{proposition}
 
The proofs are
given in section \ref{sect5} below.
An important point
is that the Fibonacci numbers
grow exponentially fast and so, for any constant $C$, we have
$$\sum_{n\geq 1}\frac{\log CF_n^2}{F_n}<\infty.$$
Thus, proposition~\ref{goodapprox2} tells us that the contribution of ``bad approximants'' to the sum
defining $B(\a)$ is universally bounded.

\remark Here, it is not critical to have an optimal bound in
proposition~\ref{constantC}. Having $c/q^2$ instead of $1/q^3$
would only replace $\sum \frac{\log F_n^2}{F_n}$ by
$\sum \frac{\log (c F_n)}{F_n}$.
\vskip.2cm

\begin{definition}
For each $\a\in B_i$, set
$$S_i(\a)=\left\{\delta\in
  U_i~\Big|~\frac{p_{n_i}}{q_{n_i}}+\delta^{q_{n_i}}=\a\right\}
\quad{\rm and}\quad
{\cal C}_i(\a) = \chi_{p_{n_i}/q_{n_i}}(S_i(\a)).$$
Moreover, define by induction 
$$V_0(\a) = \C\setminus {\cal R}_0(\a)
\quad{\rm and}\quad
V_i(\a) = V_{i-1}(\a)\setminus {\cal C}_i(\a)$$
and set 
$$S_i=S_i(\a_0),\quad {\cal C}_i ={\cal C}_i(\a_0)\quad{\rm and}\quad
V_i =V_i(\a_0).$$
\end{definition}

For $i\geq 1$, the sets ${\cal C}_i(\a)$ form periodic cycles for
$P_{\a}$ of period $q_{n_i}\geq 2$. They are clearly contained in
$V_0(\a)$ since no periodic cycle which is not fixed can belong to the
closure of the fixed external ray ${\cal R}_0(\a)$.

\vskip.2cm
\noindent{\bf Step 4.}
The conformal radius of $V_0$ may be estimated as follows.

\begin{proposition}\label{log8pi}
We have
$$\log\rad(V_0) <  -\frac{\log q_1}{q_0}+\log(8\pi).$$  
\end{proposition}

\proof 
Since $V_0(\a)$ is simply connected and avoids the fixed point 
$1-e^{2i\pi \a}$, Koebe one-quarter theorem yields
$$\log\rad(V_0) \leq \log|1-e^{2i\pi\a_0}|+\log 4.$$
We have $q_0=1$ and since $q_1=\lfloor 1/\a_0\rfloor$, we have
$$\log \left|1-e^{2i\pi\a_0}\right|<\log (2\pi \a_0) <-\frac{\log
  q_1}{q_0}+\log (2\pi).$$
\qedprop

\vskip.2cm
\noindent{\bf Step 5.}
Our goal is then to show the following inequality.
It is the main estimate of the article. 

\begin{proposition}\label{mainestimate}
For $i\geq 1$, we have
$$\log \frac{\rad(V_i)}{\rad(V_{i-1})} \leq -\frac{\log
  q_{1+n_i}}{q_{n_i}} + \frac{\log 24 F_{n_i}^2}{F_{n_i}} +
  \frac{\log 16}{1 + F_{n_i}/1.5}.$$
\end{proposition}

Combining the results from propositions \ref{goodapprox2}, \ref{log8pi}
and \ref{mainestimate}, we get (using $2<24$)
\begin{eqnarray*}
\log \rad(V_N) & < & \log \rad(V_0) +
\sum_{\underset{\scriptstyle n_i\leq N}{i\geq 1}} 
\left( -\frac{\log q_{1+n_i}}{q_{n_i}}
 +\frac{\log 24 F_{n_i}^2}{F_{n_i}} + \frac{\log 16}{1 + F_{n_i}/1.5}\right) \\
& < & -\sum_{n=0}^N \frac{\log q_{n+1}}{q_n}+\log(8\pi)+
\sum_{n=1}^N \left( \frac{\log 24 F_n^2}{F_n} + \frac{\log 16}{1 +
    F_{n}/1.5}\right) \\
& < & -\sum_{n=0}^N \frac{\log q_{n+1}}{q_n}+ 16.
\end{eqnarray*}
The proof is then completed.

Let us now explain how we get proposition \ref{mainestimate}.
First, ${\cal C}_i\subset V_{i-1}$ is the image of $S_i$ by the holomorphic
function $\chi_{p_{n_i}/q_{n_i}}$. Proposition \ref{mainestimate} is therefore
almost a consequence of the following two propositions. 

\begin{proposition}\label{radestimate}
We have the inequality
$$\log\frac{\rad(U_i\setminus S_i)}{\rad(U_i)} \leq -\frac{\log
  q_{1+n_i}}{q_{n_i}} + \frac{\log 24 F_{n_i}^2}{F_{n_i}}.$$
\end{proposition}

This proposition only consists in estimating the conformal radius of the unit
disk minus $q$ points equidistributed on a circle of radius $r<1$. The proof
is given in section \ref{sect6} below. 

\begin{proposition}\label{raddecay}
Assume $U,V\subset \C$ are two hyperbolic domains
containing $0$ and $\chi:U\to
V$ is a holomorphic map fixing $0$. Let $S$ be a finite
subset of $U$ avoiding $0$, such that $\chi(S)$ avoids $0$.
Then, 
$$\frac{\rad(V\setminus \chi(S))}{\rad(V)}\leq \frac{\rad(U\setminus
  S)}{\rad(U)}.$$
\end{proposition}

The proof of this inequality is a refinement of Schwarz's lemma and is
based on the use of ultrahyperbolic metrics (see section \ref{sect7}
below).

Combining those two inequalities would yield proposition~\ref{mainestimate} if
$\chi_{p_{n_i}/q_{n_i}}:U_i\to \C$ took its values in $V_{i-1}$, which is
almost the case.
In fact, let $\a(\delta) = p_{n_i}/q_{n_i}+\delta^{q_{n_i}}$.
As $\delta$ varies in $U_i$,
$\chi_{p_{n_i}/q_{n_i}}(\delta)$ belongs
to $V_{i-1}(a(\delta))$, which depends on $\delta$.
In section~\ref{sect8}, using Slodkowski's theorem and the
straightening of Beltrami forms, we define for $\a\in D_i$
an analytic family of universal coverings
$\pi_\a : \wt{V}_\a \to V_{i-1}(\a)$, where $\wt{V}_\a$ are open
subsets of $B(0,4)$, and $\wt{V}_{\a_0} = \D$.
The map $\chi_{p_{n_i}/q_{n_i}}(\delta)$ ``lifts'' to a map 
$\wh{\phi}(\delta)$ such that $\wh{\phi}(\delta) \in \wt{V}_{\a(\delta)}$.
It follows from the definitions that,
$$\log \frac{\rad(V_i)}{\rad(V_{i-1})} = \log
\frac{\rad(\wt{V}_{\a_0} \setminus \pi_{\a_0}^{-1}(\phi(S)))
}{\rad(\wt{V}_{\a_0})}.$$
Now $\wt{V}_{\a_0} = \D$ and $\wh{\phi}(S) \subset
\pi_{\a_0}^{-1}(\phi(S))$, thus
$$\log \frac{\rad(V_i)}{\rad(V_{i-1})} \leq \log
\rad(\D \setminus \wh{\phi}(S)).$$
We would like to apply proposition~\ref{raddecay} to $U=B_i$ and $\chi =
\wh{\phi}$. But we cannot take $V=\D$ because $\wh{\phi}$ does not necessarily
take its values in $\D$. 
However, in section~\ref{sect8}, we prove the following estimate.
\begin{proposition}\label{movingrange}
For $\a \in B_i$, the sets $\wt{V}_\a$ are all
contained in some ball $B(0,\rho_2)$ with
$$\log \rho_2= \frac{\log 16}{1+q_{n_i}/1.5}$$
\end{proposition}

We can therefore take $V=B(0,\rho_2)$ in
proposition~10 and we obtain

$$\log  \frac{\rad(B(0,\rho_2)\setminus \wh\phi(S))}{\rad(B(0,\rho_2))}
\leq \log\frac{\rad(U_i\setminus S)}{\rad (U_i)}.$$
Now, by inclusion
\begin{eqnarray*}
\log \rad(\D\setminus \wh\phi(S)) & \leq &
\log \rad(B(0,\rho_2)\setminus \wh\phi(S)) \\
& \leq & \log\frac{\rad(U_i\setminus S)}{\rad (U_i)} +
\log \rad(B(0,\rho_2))\\
& \leq & -\frac{\log q_{1+n_i}}{q_{n_i}} + \frac{\log 24
F_{n_i}^2}{F_{n_i}}
+\frac{\log 16}{1+q_{n_i}/1.5}
\end{eqnarray*}
according to propositions ~\ref{radestimate} and~\ref{movingrange}.

\section{Parabolic Explosion.\label{sect3}}

The proofs of propositions \ref{functionchi} and \ref{constantC} may be 
found in \cite{c} or \cite{bc},
but for completeness, we sketch them here. 

\vskip.2cm
\noindent{\bf Proof of Proposition \ref{functionchi}.}
It is well known that, when $\a$ varies,
periodic points with multiplier different from $1$
can be locally followed holomorphically in terms of $\a$.
To prove this, one applies the implicit function theorem
(complex-analytic version) to the equation
``$P_{\a}^{\circ k}(z)-z = 0$''
where $k$ is the period: at a point $(\alpha,z)$ in the surface
defined by, the derivative with respect to $z$ is equal to $m-1$ where $m$
is the multiplier. In our case, on
$B'=B(\frac{p}{q},R(\frac{p}{q}))\setminus\{\frac{p}{q}\}$,
no points of period dividing $q$ has multiplier $1$. 

Since $B'$ is not simply connected, the holomorphic
dependence in terms of $\a$ may have a monodromy when $\a$ makes one turn
around $p/q$.
Let $B=B(\frac{p}{q},R(\frac{p}{q}))$. Let us consider the subset of
$B\times \C$ defined by 
\[\cal M = \big\{(\a,z)\big| P_\a^{\circ q}(z)-z=0 \big\}.\]
For $\a$ at the center of $B$ (i.e.\ $\a=p/q$),
only $z=0$ is parabolic. Thus
the fixed points of $P_{p/q}^{\circ q}$ that are different from $0$
have no monodromy: they can be followed holomorphically as a function of $\a$
on all of $B$. The graphs of these functions
$\a \mapsto z(\a)$ are connected components of $\cal M$.
There is only one other component. It contains $(\a,z)=(p/q,0)$, and it is
singular. To study it, one looks at the expansion of the equation at
this point.
First, it is known (see \cite{dh},
chapter IX) that there exists a complex number $A\in \C^*$ such that 
\[P_{p/q}^{\circ q}(z) = z+Az^{q+1} + {\cal O}(z^{q+2}).\]
This means there are $q+1$ fixed points of $P_{p/q}^{\circ q}$ at $z=0$.
Then,
\[P_{p/q+\e}^{\circ q}(z)-z = z\cdot\big(2i\pi q \e +
A z^q + O(z\e)+O(z^{q+1})\big).\]
To get rid of the singularity, one considers a new variable 
$\delta\in D=B(0,R(p/q)^{1/q})$
related to $\a$ by $\a = \delta^q + p/q$. This transforms the component
of $\cal M$ containing 
$(p/q,0)$ into the union of $1+q$ graphs of functions from 
$D$ to $\C$ that are transversal and meet only at $(0,0)$. This can be
proved by blowing-up $D \times \C$ at $(0,0)$, 
i.e. by introducing a new variable, the slope $\lambda = z/d$.
One of the graphs corresponds to the fixed point $z=0$ of $P_{\a}$
which does not move, and the others are graphs of functions 
$\phi_1(\delta)$, $\phi_2(\delta)$, \ldots, $\phi_q(\delta)$ passing through
$(0,0)$ with slopes $\lambda$ equal to the $q$-th roots of $-2i\pi
q/A$. 

The function $\chi$ of Proposition \ref{functionchi} is any of these functions
$\phi_i$. The points $P_\alpha(\chi(\delta))$ and $\chi(\zeta \delta)$
are both fixed points of $P_\alpha^{\circ q}$. Since there graphs pass
through $(0,0)$ they coincide with functions $\phi_i$ and
$\phi_j$ for some $i$ and $j$. By comparing the derivatives at $\delta=0$, one
gets $i=j$ and so
$$ P_\alpha(\chi(\delta)) = \chi(\zeta \delta).$$
This shows that the set
$$\Big< \chi(\delta),\chi(\zeta\delta), 
\chi(\zeta^2\delta),\ldots,\chi(\zeta^{q-1}\delta)\Big>$$
forms a cycle of period $q$ of $P_{p/q+\delta^q}$.
\qedprop

\vskip.2cm
\noindent{\bf Proof of proposition \ref{constantC}.}
The only values of the parameter $\a$ for which $P_\a$ has a
multiple fixed point are the integers. Thus, the result is trivial for
$q=1$. Let us now assume that $q\geq 2$. 
In that case, the proof relies on Douady's landing theorem and 
the Pommerenke-Levin-Yoccoz inequality (see \cite{h} or \cite{pe}). 

\begin{figure}[htbp]
\begin{center}
 \begin{picture}(150,160)%
  \put(0,0){\scalebox{0.5}{%
   \includegraphics*[155,390][455,670]{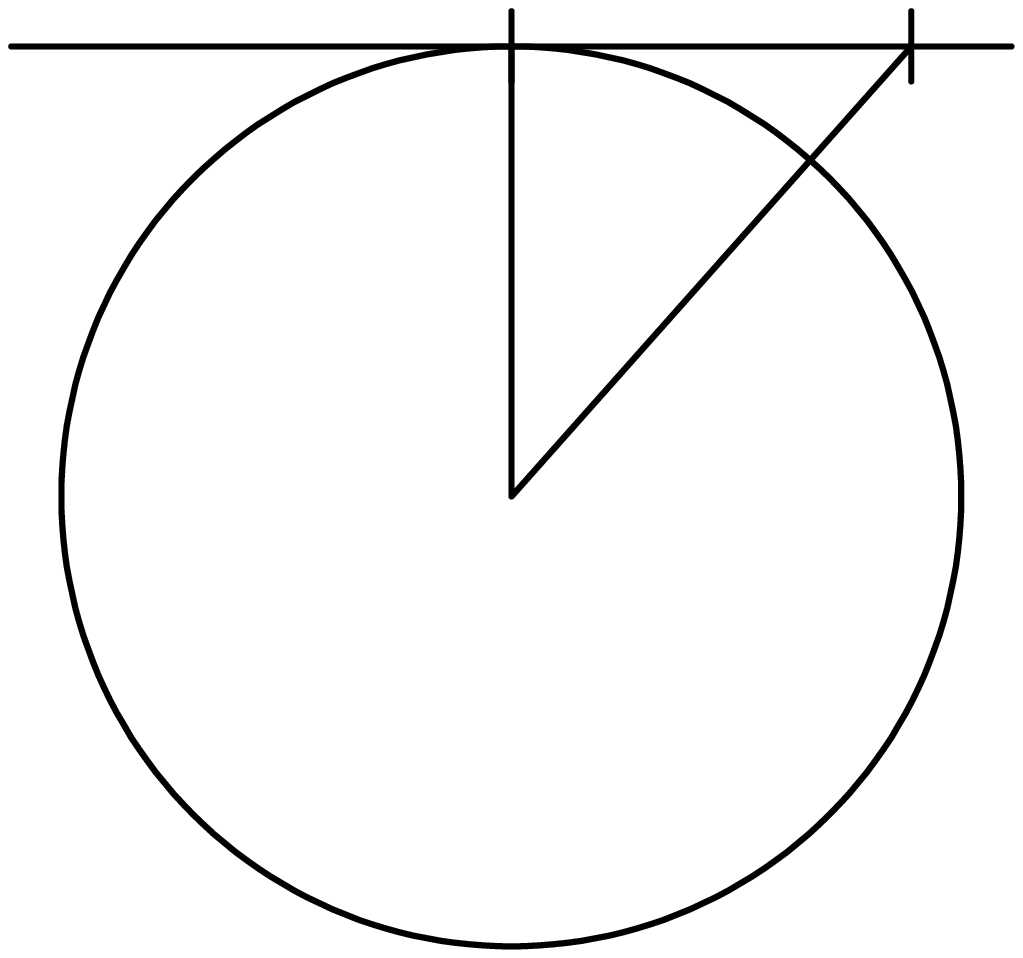}
  }}
  \put(66,148){$p'/q'$}
  \put(127,148){$p/q$}
  \put(155,130){$\R$}
  \put(-50,50){Yoccoz disk}
 \end{picture}
\end{center}
\caption{The complex number $\a$ lies somewhere in the Yoccoz disk}
\end{figure}

Let us choose a rational number $p/q$, and assume that
$\a\neq p/q$ and $P_\a^{\circ q}$ has a multiple fixed point
$z_0$. Then, $P_\a$ has a parabolic cycle
$\left<z_0,z_1,\ldots, z_{q_1-1}\right>$ of period $q_1$ dividing $q$,
and the immediate basin of this parabolic cycle contains
the critical point $\omega_0=-e^{2i\pi\a}/2$ of $P_\a$. As a
consequence, the Julia set $J(P_\a)$ is connected and all other periodic 
cycles of $P_\a$ are repelling.

If $0$ is parabolic, then $\a=p'/q$
with $p'$ not necessarily prime to $q$. So, the distance between $\a$ and
$p/q$ is bounded from below by $1/q$.

Otherwise, $0$ must be repelling, thus $\a$ belongs
to the lower half-plane $\{\im(\a)<0\}.$  Since the Julia set is
connected, Douady's landing theorem asserts that there are finitely
many rays landing at $0$, let's say $q'$. Those $q'$ rays can be
ordered cyclically by their arguments. They are permuted by $P_\a$
and each ray is mapped to the one which is $p'$ further
counter-clockwise for some $p'<q'$, $p'$ prime to $q'$.
Then, the Yoccoz inequality
implies that $\a$ belongs to the closed disk of radius $(\log2)/(2\pi
q')$ tangent to the real axis at $p'/q'$ (see for example \cite{h}).
A key combinatorial lemma that is proved below is that we necessarily
have $q>q'$.

The Pythagoras theorem then gives
$$\left|\a-\frac pq\right|\geq \sqrt{
\left(\frac{p'}{q'}-\frac{p}{q}\right)^2+
\left(\frac{\log 2}{2\pi q'}\right)^2}-\frac{\log 2}{2\pi q'}.$$
Since $q>q'$, and $q\geq 2$, an elementary computation gives
$|\a-p/q|>1/q^3$. \qedprop

\begin{lemma}[Key combinatorial lemma]
$$q' < q$$
\end{lemma}

\begin{figure}[htbp]
\begin{center}
 \begin{picture}(225,175)%
  \put(0,0){\scalebox{0.5}{%
   \includegraphics*[80,220][530,570]{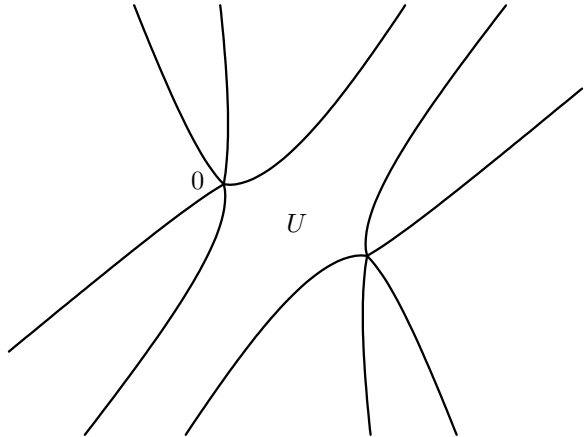} 
  }}
  \put(110,84){$U$}
  \put(74,100){$0$}
 \end{picture}
\end{center}
\caption{Schematic example for $q'=5$.}
\end{figure}

\proof By assumption, $q>1$.
Consider the 
complement in $\C$ of the $q'$ external rays landing at $0$ together
with this point. It has $q'$ connected components.
Let $V$ be the one containing the critical
point. The orbit of the critical point of $P_\a$ must first
visit each
connected component of this complement, before first falling back
somewhere in $V$. Since the critical point belongs to the immediate
basin of the parabolic cycle, this implies the period is $\geq q'$,
and thus $q \geq q'$.

Let us assume by contradiction that $q=q'$. The point $0$ has two
distinct preimages: $0$ and another point. Consider the union of these
two points and the $2q$ external rays landing at them. 
Let $U$ be the component of the complement of this union containing the critical point. 
It is known that 
$P_\a^{\circ q}(U)=V$ and $P_\a^{\circ q} : U \to V$ is a proper
ramified covering of degree $2$. Let $f$ be the
restriction $P_\a^{\circ q} : \overline{U} \to \overline{V}$. 
Note that $\overline{U} \subset \overline{V}$. 
The contradiction follows from a version of
the Lefschetz fixed point formula (see lemma~3.7 in \cite{gm}): the
point $z=0$ is 
fixed, and the point of the parabolic cycle whose immediate basin
contains the critical point is a multiple fixed point of
$f$. Thus the sum of Lefschetz indices is $\geq 3$, whereas
according to the Lefschetz formula is should be equal to the degree of
$f$, i.e., $2$. This leads to contradiction and thus $q>q'$.
\qedprop

\section{proof of proposition \ref{constantC'}.\label{sect4}}

The set ${\cal B}$ is contained in the union of $\Z$ and the lower
half-plane. It is the graph of the function 
$$f:x\mapsto\frac{1}{2\pi}\log \cos(2\pi x)$$
which is periodic of period $1$ and defined for $x\in ]-1/4+k,1/4+k[$,
$k\in \Z$. We will now show that for $x\in ]0,1/2[$,
$f(x-x^2)<-x^2$. Since $f$ 
is decreasing on $[0,1/2[$, proposition \ref{constantC'} follows. 
\begin{figure}[htbp]
\begin{center}
 \begin{picture}(222,170)%
  \put(0,0){\scalebox{0.5}{%
   \includegraphics*[80,130][524,470]{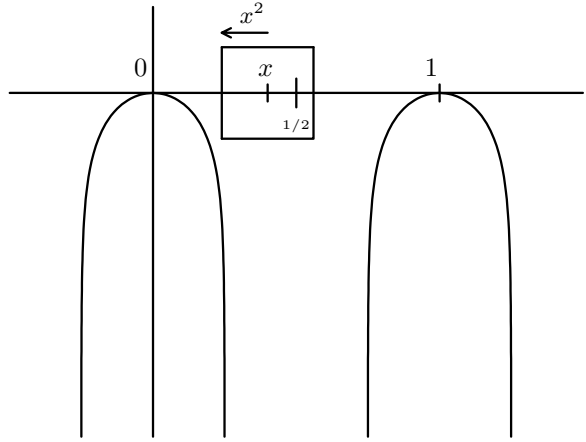}
  }}
  \put(99,140){$x$}
  \put(92,160){\small $x^2$}
  \put(52,140){$0$}
  \put(108,120){\tiny $1/2$}
  \put(162,140){$1$}
 \end{picture}
\end{center}
\caption{The graph of the function $f$. The inequality $f(x-x^2)<-x^2$
  means the lower left corner of the square is above the curve.}
\end{figure}
We want to show that the function $g(x) = x^2+f(x-x^2)$ is negative on
$]0,1/2[$. Since $g(0)=0$, it is sufficient to show that $g'(x)<0$ on
$]0,1/2[$. This is equivalent to proving that 
$$\tan(2\pi(x-x^2))>\frac{2x}{1-2x}.$$
Let us make the change of variable $u=1/2-x$. Then, the previous becomes
$$\tan(2\pi u^2)<\frac{2u}{1-2u}.$$
We are done since for all $u\in ]0,1/2[$, we have
$$\sin (2\pi u^2) < 2\pi u^2\quad {\rm and}\quad
\cos(2\pi u^2) > \frac{\pi/2-2\pi u^2}{\pi/2-0}=1-4u^2 >\pi u(1-2u).$$
\qedprop 

\section{Good approximants.\label{sect5}}

Note that we made the assumption $\a_0\in ]0,1/2[$ and so $q_1\geq 2$. In
particular, for all $i\geq 1$, we have $q_{n_i}\geq 2$. 

The inclusion $D_1\subset \C\setminus {\cal B}$ follows from proposition 
\ref{constantC'}.
For $i\geq 1$, the inclusion $B_{i}\subset D_{i}$ is immediate since the two
disks have the same center and the radius of $D_{i}$ is $q_{n_i}$ times the
radius of $B_{i}$. 

The classical estimate we will use is that for all $n\geq 0$, we have
$$\frac{1}{2q_nq_{n+1}} < \left|\a_0-\frac{p_n}{q_n}\right|<
\frac{1}{q_nq_{n+1}}.$$
Thus, we have 
$$\left|\a_0-\frac{p_{n_i}}{q_{n_{i}}}\right| <
\frac{1}{q_{n_{i}}q_{1+n_{i}}}<\frac{1}{2q_{n_i}^3}$$
by definition of ${\cal N}$. In particular $\a_0$ belongs to $B_i$
and is closer to the center than to the boundary.

Moreover, for $i\geq 1$, $q_{n_{i+1}}\geq q_{1+n_i}>2q_{n_i}^2$ and
$q_{1+n_{i+1}}>2q_{n_{i+1}}^2\geq 2q_{1+n_i}^2$. So
$$\left|\a_0-\frac{p_{n_{i+1}}}{q_{n_{i+1}}}\right| < 
\frac{1}{q_{n_{i+1}}q_{1+n_{i+1}}} 
<\frac{1}{4q_{n_{i}}^2q_{1+n_i}^2}
<\frac{1}{4q_{n_i}q_{1+n_i}}
<\frac{1}{2}\left|\a_0-\frac{p_{n_i}}{q_{n_{i}}}\right|.$$
In other words, the distance from $\a_0$ to the center of
$B_{i+1}$ is less than half the distance from $\a_0$ to the center of $B_i$.
It follows from these two claims that in $B_i^*$, one can fit a
disk centered at 
$p_{n_{i+1}}/q_{n_{i+1}}$ with radius at least equal to 
$$\frac{1}{2}\left|\a_0-\frac{p_{n_i}}{q_{n_{i}}}\right|
>\frac{1}{4q_{n_i}q_{1+n_i}}> \frac{1}{q_{n_{i+1}}^2}.$$
Indeed, $q_{n_{i+1}}\geq q_{1+n_i}>2q_{n_i}^2\geq 4q_{n_i}$. 
In particular, for all $i\geq 1$, $D_{i+1}\subset B_i^*$. 
We proved proposition~\ref{goodapprox1}.

Finally, if $n\geq 1$ and $n\notin{\cal N}$, we have $q_{n+1}\leq 2q_n^2$.
It follows that 
$$\sum_{n=1}^{N} \frac{\log q_{n+1}}{q_{n}} <
\sum_{\underset{\scriptstyle n_i\leq N}{i\geq 1}} 
\frac{\log q_{1+n_i}}{q_{n_i}}
+\sum_{n\in [1,N]\setminus {\cal N}} \frac{\log 2q_n^2}{q_n}.$$
Since $q_n$ is bounded from below by the $n$-th Fibonacci number
$F_n$, according to the lemma~\ref{dulllemma} in the appendix,
we have 
$$\sum_{n\in [1,N]\setminus {\cal N}} \frac{\log 2q_n^2}{q_n}\leq
\sum_{n\in [1,N]\setminus {\cal N}} \frac{\log 2F_n^2}{F_n}.$$
We proved proposition~\ref{goodapprox2}.

\section{An estimate for a conformal radius.\label{sect6}}

Here, we prove proposition~\ref{radestimate}.

\begin{definition}
Given an integer $q\geq 1$, set
$$\U_q = \left\{e^{2i\pi k/q}~\big|~k=0,\ldots,q-1\right\}.$$
\end{definition}

The following estimate was explained to us by Douady. 

\begin{proposition}
There exists a constant $C>0$ such that for $q\geq 2$ and $r<1$, we have 
$$\log\rad(\D\setminus r\U_q) \leq \log r + \frac{C}{q}.$$
one can take $C=\log 4 + 2\log(1+\sqrt{2})$. 
\end{proposition}

\begin{figure}[htbp]
\begin{center}
 \begin{picture}(350,150)%
  \put(0,0){\scalebox{0.35}{%
   \includegraphics*[100,190][510,600]{dessin4.eps} 
  }}
  \put(200,0){\scalebox{0.35}{%
   \includegraphics*[100,190][510,600]{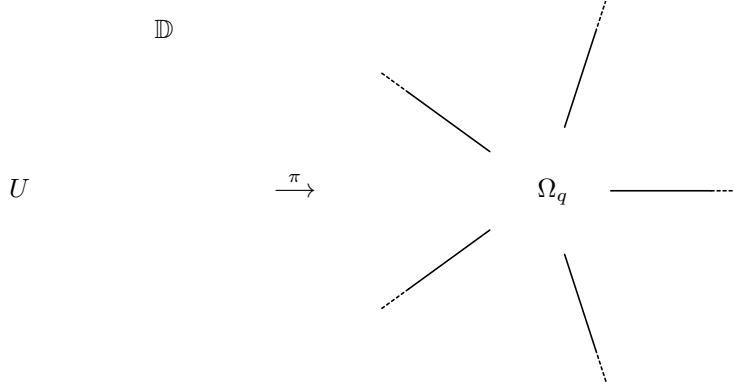} 
  }}
  \put(70,70){$U$}
  \put(125,130){$\D$}
  \put(270,70){$\Omega_q$}
  \put(170,70){$\overset{\pi}{\longrightarrow}$}
 \end{picture}
\end{center}
\caption{The map $\pi$ maps $U$ to the slit plane $\Omega_q$}
\end{figure}

\proof
By inclusion, we have 
$$\rad(\D\setminus r\U_q)\leq \rad(\C\setminus r\U_q) = r\cdot
\rad(\C\setminus \U_q).$$
Let $\pi:\D\to \C\setminus\U_q$ be a universal covering which sends $0$ to
$0$. By symmetry, for $k=0,\ldots,q-1$, the half lines 
$L_k=\{\rho e^{2i\pi k/q}~|~\rho>1\}$ are geodesics in $\C\setminus \U_q$ for
the hyperbolic metric. Set $\Omega_q = \C\setminus \bigcup L_k$. 
There is a formula for the conformal representation
$\phi_q:\D\to\Omega_q$:
$$\phi_q(z) = z\left(\frac{4}{(1-z^q)^2}\right)^{1/q}.$$
In particular, we have 
$$\rad(\Omega_q) = 4^{1/q}.$$

Now, the connected component $U$ of
$\pi^{-1}(\Omega_q)$ which contains $0$ is bounded by $2q$
geodesic arcs of circles in $\D$ whose endpoints are equidistributed on $S^1$.
An elementary computation shows that $U$ contains the disk centered at $0$
with radius 
$$\rho_q = \frac{1-\tan(\pi/4q)}{1+\tan(\pi/4q)}$$
Since the image by $\pi$ of this disk is contained in $\Omega_q$, it follows
from Schwarz's lemma that 
$$ \rad(\C\setminus \U_q) \leq \frac{\rad(\Omega_q)}{\rho_q}.$$
thus $$\log \rad(\C\setminus \U_q) \leq \frac{\log 4}{q}
+ \log \frac{1+\tan(\pi/4q)}{1-\tan(\pi/4q)}.$$
By convexity of $\ds f(x)= \log \frac{1+\tan x}{1-\tan x}$ on
$[0,\pi/8]$, we have $\ds f(\pi/4q) \leq \frac{2}{q} f(\pi/8)$.
The result now follows easily. 
\qedprop

We can now estimate the conformal radius of $U_i\setminus S_i$ for $i\geq 1$.
The radius of the ball $U_i$ is $(1/q_{n_i}^3)^{1/q_{n_i}}$ and the set
$S_i$ consists of $q_{n_i}$ points equidistributed on a circle of radius 
$$\left|\a_0-\frac{p_{n_i}}{q_{n_i}}\right|^{1/{q_{n_i}}}<
\left(\frac{1}{q_{n_i}q_{1+n_{i}}}\right)^{1/{q_{n_i}}}.$$
So, we have
\begin{eqnarray*}
\log\frac{\rad(U_i\setminus S_i)}{\rad(U_i)} & < & 
\log \frac{(1/q_{n_i}q_{1+n_{i}})^{1/{q_{n_i}}}}
{(1/q_{n_i}^3)^{1/q_{n_i}}} +
\frac{C}{q_{n_i}}\\ 
 & = & - \frac{\log q_{1+n_i}}{q_{n_i}} + 2 \frac{\log{q_{n_i}}}{q_{n_i}} +
 \frac{C}{q_{n_i}}.
\end{eqnarray*}
Since $q_n$ is bounded from below by the $n$-th Fibonacci number $F_n$, we
have
$$2 \frac{\log{q_{n_i}}}{q_{n_i}} + \frac{C}{q_{n_i}}
\leq \frac{\log{24 q_{n_i}^2}}{q_{n_i}} \leq
\frac{\log{24 F_{n_i}^2}}{F_{n_i}}$$
according to lemma~\ref{dulllemma} in the appendix.

\section{Comparison between conformal radii.\label{sect7}}

Our goal in this section is to prove proposition \ref{raddecay}.
The proof relies on a relative Schwarz's lemma. 

\subsection{A relative Schwarz's lemma.}

\begin{definition}
A metric $\rho|dz|$, $\rho\geq 0$ is said to be ultrahyperbolic in a
Riemann surface $X$ if it has the following properties:
\begin{itemize}
\item[(i)] $\rho$ is upper semicontinuous.
\item[(ii)] At every $x_0\in X$ with $\rho(x_0)>0$ there exists a
``supporting metric'' $\rho_0$, defined and of class $C^2$ in a
neighborhood $V$ of $x_0$, such that $\Delta \log \rho_0\geq \rho_0^2$
and $\rho\geq \rho_0$ in $V$, while $\rho$ coincides with $\rho_0$ at $x_0$. 
\end{itemize}
\end{definition}

In a hyperbolic Riemann surface $X$,
there exists a unique maximal ultrahyperbolic metric $\rho_X$, and
this metric has constant curvature $-1$. 
It is maximal in the sense that every ultrahyperbolic metric $\rho$ on $X$
satisfies $\rho\leq \rho_X$ throughout $X$. 
This maximal metric is called the \em Poincar\'e metric \em on $X$. 

For example, the Poincar\'e metric $\rho_\D$ on the unit disk $\D$ is 
$$\rho_\D=\frac{2}{1-|z|^2}|dz|.$$
More generally, if $\pi:\D\to X$ is a universal covering, the
Poincar\'e metric $\rho_X$ coincides with the unique metric such that
$\pi^*\rho_X=\rho_\D$. 

Now, if $f:\D\to \D$ is holomorphic, then $f^*\rho_\D$ is
ultrahyperbolic on $\D$, and thus $f^*\rho_\D\leq \rho_\D$. This may
be written as the Schwarz-Pick theorem:
$$(\forall z\in \D)\quad |f'(z)|\leq \frac{1-|f(z)|^2}{1-|z|^2}.$$
More generally, if $X$ and $Y$ are two hyperbolic Riemann surfaces,
then every holomorphic map $f:(X,\rho_X)\to (Y,\rho_Y)$ is
contracting:
$$f^*\rho_Y \leq \rho_X.$$
In particular, if $X\subset Y$, then $\rho_Y\leq \rho_X$.

In this subsection, we are interested in comparing the relative
contraction of a holomorphic map $f:X\to Y$ for several Poincar\'e
metrics. We will show that if $Y'\subset Y$ is an arbitrary open
subset and $X'=f^{-1}(Y')$, then $f:(X',\rho_{X'})\to (Y',\rho_{Y'})$
is less contracting than $f:(X,\rho_X)\to (Y,\rho_Y)$. 

\begin{lemma}{\sc (Relative Schwarz's Lemma)}
Let $f:X\to Y$ be an analytic map between
two hyperbolic Riemann surfaces. Let $Y'\subset Y$ be an arbitrary
open subset and set $X'=f^{-1}(Y')$. Then, on $X'$,
$$\frac{f^*\rho_Y}{\rho_X}\leq \frac{f^*\rho_{Y'}}{\rho_{X'}}\leq 1.$$
\end{lemma}

The main tool of the proof is the use of Ahlfors's ultrahyperbolic
metrics (see\cite{a} for example). This was suggested to us by
McMullen. 

\proof
Let us first consider the case where $y_0\in Y$ is an arbitrary point
and $Y'=Y\setminus\{y_0\}$. We will show that the metric $\sigma$
defined on $X'=X\setminus f^{-1}\{y_0\}$ by 
$$\sigma = \frac{f^*\rho_Y}{f^*\rho_{Y'}}\rho_{X'}$$
extends continuously to a ultrahyperbolic metric on $X$. It  will then follow
from the definition of the Poincar\'e metric $\rho_X$ that 
$$\frac{f^*\rho_Y}{f^*\rho_{Y'}}\rho_{X'}\leq \rho_X,$$
which is the required result. 

\vskip.2cm
\noindent{\bf Step 1.}
The metric $\sigma$ is a priori only defined on $X'\setminus {\rm
Crit}(f)$, where ${\rm Crit}(f)$ is the set of critical points of
$f$. But since
$$\sigma=\left(\frac{\rho_Y}{\rho_{Y'}}\circ f\right)\cdot
\rho_{X'},$$
we see that $\sigma$ is positive and $C^2$ on $X'$. We will now show
that $\Delta\log \sigma\geq \sigma^2$ on $X'\setminus {\rm
Crit}(f)$. Since ${\rm Crit}(f)$ is discrete in $X'$, this inequality
holds on $X'$. Therefore, $\sigma$ is ultrahyperbolic on $X'$. 

On $X'\setminus {\rm Crit}(f)$ we have
$$\Delta \log \sigma = \Delta \log f^*\rho_Y + \Delta\log \rho_{X'}
-\Delta\log f^*\rho_{Y'}
= [f^*\rho_Y]^2+[\rho_{X'}]^2-[f^*\rho_{Y'}]^2.$$
The second equality comes from the fact that the three metrics have
curvature $-1$. 
Since $Y'\subset Y$, we have $\rho_Y\leq \rho_{Y'}$ on $Y'$, and so,
$f^*\rho_Y\leq f^*\rho_{Y'}$ on $X'$. And since, $f^*\rho_{Y'}$ is
ultrahyperbolic on $X'$, we have $f^*\rho_{Y'}\leq \rho_{X'}$ on
$X'$. Now, if $a$, $b$ and $c$ are three positive numbers such that
$a\leq b\leq c$, then
$$(c-b)(b-a) \geq 0 \quad \Longrightarrow \quad
cb - b^2+ba \geq ac \quad \Longrightarrow \quad 
a+c-b\geq \frac{ac}{b}.$$
Therefore, 
$$\Delta \log \sigma =[f^*\rho_Y]^2+[\rho_{X'}]^2-[f^*\rho_{Y'}]^2\geq 
\left[\frac{f^*\rho_Y\cdot \rho_{X'}}{f^*\rho_{Y'}}\right]^2=\sigma^2.$$

\vskip.2cm
\noindent{\bf Step 2.}
We claim that we may extend $\sigma$ continuously to $X\setminus X'$
by setting $\sigma = f^*\rho_Y$ there. Indeed, let $x_0$ be an
arbitrary point in $X\setminus X'$. It is sufficient to show that
$$\lim_{x\to x_0}\frac{\rho_{X'}}{f^*\rho_{Y'}}(x) = 1.$$

\begin{lemma}
 Let $X$ be a hyperbolic Riemann surface, $X'$ be an open subset
 of $X$ and assume that $x_0 \in X\setminus X'$ is an isolated point
 of $X \setminus X'$.
 Then, in any analytic chart, if we note
 $\rho_{X'} = \rho_{X'}(x) |dx|$, we have
 $$\rho_{X'}(x) \underset{x\to x_0}\sim \frac{1}{|x-x_0|\log\ds\frac{1}{|x-x_0|}}.$$
\end{lemma}
Note that the formula is independent of the chosen chart.

\proof
We will work in the local coordinates given by the universal covering
$\pi_X:(\D,0)\to (X,x_0)$. We set $U=\pi_X^{-1}(X')$. Let $x=\pi_X(w)$
Then, we have $\rho_{X'} = \rho_U(w)|dw|$. We claim that
$$\rho_U(w)\underset{w\to 0}\sim
\frac{1}{|w|\log\ds\frac{1}{|w|}}.$$ 
Indeed, we may find $\e>0$ such that $\D_\e^*\subset U\subset
\D^*$. Then, 
$$\rho_{\D^*}(w) = \frac{1}{|w|\log\ds\frac{1}{|w|}}
\leq \rho_U(w) \leq \rho_{\D_\e^*}(w)
= \frac{1}{|w|\log\ds\frac{\e}{|w|}}.$$
Then, since $\rho_{X'}(x) = \rho_{X'}(\pi_X(w)) = \rho_U(w) / |\pi'_X(w)|$, and
$|x-x_0|\sim |\pi'_X(0)|\cdot |w|$ when $w \to 0$, the result follows.
\qedlem

Now, let us choose analytic charts for $X$ near $x_0$
and for $Y$ near $y_0$, and note $r=|x-x_0|$.
Then, as $x \to x_0$, we have 
$$|f(x)-f(x_0)|  \sim  Ar^d \quad{\rm and}\quad
|f'(x)| \sim Adr^{d-1},$$
where $d$ is the local degree of $f$ at $x_0$ and $A>0$.
Thus, as $x \to x_0$, we have
$$\frac{\rho_{X'}}{|dx|}(x) \sim \frac{1}{r \log \ds\frac{1}{r}},$$
and 
$$\frac{f^*\rho_{Y'}}{|dx|}(x) = |f'(x)| \rho_Y(f(x))
\sim \frac{A d r^{d-1}}{A r^d \log \ds\frac{1}{A r^d}} \sim \frac{1}{r \log \ds\frac{1}{r}}.$$
The claim follows.

\vskip.2cm
\noindent{\bf Step 3.}
As we have just seen, we may extend $\sigma$ continuously to
$X\setminus X'$ by setting $\sigma = f^*\rho_Y$ there, and since
$f^*\rho_{Y'}\leq \rho_{X'}$, we see that $\sigma\geq f^*\rho_Y$. If
$\sigma$ does not vanish at $x_0$, i.e., if $x_0$ is not a critical
point of $f$, then $f^*\rho_Y$ is $C^2$ in a neighborhood of $x_0$, has
curvature $-1$ and coincides with $\sigma$ at $x_0$. Thus, condition
(ii) in the definition of ultrahyperbolic metrics is satisfied:
$f^*\rho_Y$ is a ``supporting metric'' at $x_0$ and we have proved
that $\sigma$ is ultrahyperbolic. Thus, if $Y'$ is obtained by
removing one point from $Y$, the relative Schwarz's lemma is proved. 

By induction, the lemma is therefore proved when $Y'$ is obtained by
removing finitely many points from $Y$. In order to prove the lemma
for an arbitrary open subset $Y'\subset Y$, we may choose a dense
countable set $\{y_n,~n\geq 0\}\subset Y\setminus Y'$, define
$Y_n=Y\setminus \{y_k,~k\leq n\}$ and set $X_n= f^{-1}(Y_n)$. Then,
for all $n\geq 0$, 
\begin{equation}\label{relcont}
\frac{f^*\rho_Y}{\rho_X}\leq \frac{f^* \rho_{Y_n}}{\rho_{X_n}}\leq
1.
\end{equation}

\begin{lemma}
Assume $(U_n)_{n\geq 0}$ is a decreasing sequence of hyperbolic Riemann
surfaces. Let $U$ be the interior of $\ds \bigcap_{n\geq 0} U_n$. 
As $n\to +\infty$, the Poincar\'e metrics $\rho_{U_n}$ converge uniformly on
every compact subset of $U$ to the Poincar\'e metric $\rho_U$.
\end{lemma}

\proof
Let $a$ be an arbitrary point in $U$ and let $U_a$ be the connected component
of $U$ that contains $a$. Let
$\phi_n:(\D,0)\to (U_n,a)$ and $\phi:(\D,0)\to (U_a,a)$ be the universal
coverings which have real and positive derivatives at $0$ (for some chart
around $a$ in $U_a$). We will show that the maps $\phi_n$ converge 
uniformly to $\phi$ on every compact subset of $U_a$. The lemma follows
easily. 

The maps $\phi_n$ all take their values in $U_0$ which is
hyperbolic. So, they form a normal family. Let $\psi:(\D,0)\to (U_a,a)$ be a
limit value. For all $n\geq 0$, the map $\phi$ takes its values in $U_n$, and
thus, $\phi'(0)\leq \phi_n'(0)$. Similarly, the map $\psi$ takes its values in
$U_a$ and thus, $\psi'(0)\leq \phi'(0)$. Since $\psi$ is a limit value of
the sequence $\phi_n$, we have $\psi'(0) = \phi'(0)$ and by the classical
Schwarz lemma,
$\psi=\phi$.  
\qedlem

As a consequence, as $n\to +\infty$, the Poincar\'e metrics $\rho_{X_n}$ and
$\rho_{Y_n}$ converge uniformly on every compact subset of $X'$ and $Y'$ to the
Poincar\'e metrics $\rho_{X'}$ and $\rho_{Y'}$.
Passing to the limit in inequality (\ref{relcont}) gives the required result: 
 $$\frac{f^*\rho_Y}{\rho_X}\leq \frac{f^* \rho_{Y'}}{\rho_{X'}}\leq
1.$$
\qedprop

\subsection{Proof of proposition \ref{raddecay}.}
Let us recall the problem. We assume $U,V\subset \C$ are hyperbolic
domains containing $0$, we assume $\chi:(U,0)\to (V,0)$ is holomorphic, and
we assume that $\chi(S)$ avoids $0$ (in which case $S$ also avoids $0$). 
We wish to conclude that 
$$\frac{\rad(V\setminus \chi(S))}{\rad(V)}\leq 
\frac{\rad(U\setminus S)}{\rad(U)}.$$
The conformal radius $\rad(U)$ is related to the coefficient of the
Poincar\'e metric $\rho_U(0)$ as follows:
$$\rad(U) = \frac{2}{\rho_U(0)}.$$
We will apply the relative Schwarz's lemma with $X=U$, $Y=V$, $f=\chi$,
$Y'=V\setminus \chi(S)$ and $X'=\chi^{-1}(Y')$. We have
$$\frac{f^*\rho_Y}{\rho_X}\leq \frac{f^*\rho_{Y'}}{\rho_{X'}},$$
which may be rewritten as
$$\frac{\rho_{X'}}{\rho_X}\leq
\frac{f^*\rho_{Y'}}{f^*\rho_Y}=\frac{\rho_{Y'}}{\rho_Y}\circ f.$$
Evaluating this inequality at $0$, and using the relation between the
conformal radius and the coefficient of the Poincar\'e metric, we get
$$\frac{\rad(V\setminus \chi(S))}{\rad(V)} = \frac{\rad(Y')}{\rad(Y)}
\leq \frac{\rad(X')}{\rad(X)} = \frac{\rad(U\setminus
  \chi^{-1}(\chi(S)))}{\rad(U)}.$$
The result follows since $U\setminus \chi^{-1}(\chi(S))\subset U\setminus
S$, and so,
$$\rad(U\setminus \chi^{-1}(\chi(S)))\leq \rad(U\setminus S).$$

\section{Holomorphic motions.\label{sect8}}

To prove proposition~\ref{mainestimate},
we must now take into account the fact that for $i\geq
2$, $\chi_{p_{n_i}/q_{n_i}}$ does not take its values in $V_{i-1}$ 
but rather that
$\chi_{p_{n_i}/q_{n_i}}(\delta)$ belongs to 
$V_{i-1}(\a(\delta))$ with
$\a(\delta)=p_{n_i}/q_{n_i}+\delta^{q_{n_i}}$. 
The sets $V_{i-1}(\a)$ move holomorphically with respect to $\a\in
D_i$ and when $\delta$ ranges in $U_i$, $\a(\delta)$ remains
in $B_i$
which is well inside $D_i$ (the ratio of the radii is
$q_{n_i}$ and $q_{n_i} \geq 2$ as $\a \in\,]0,1/2[$).

To begin with, let us work in quite a general, but normalized, setting under
the following assumptions.  
We assume that $V_\lambda$ are hyperbolic subdomains of $\C$ which contain $0$
and move holomorphically with respect to $\lambda\in \D$.
By Slodkowski's theorem, we can assume that the holomorphic motion is a 
holomorphic motion of the whole complex plane. We set 
$${\cal V} =\big\{(\lambda,z)~\big|~\lambda\in \D~{\rm and}~z\in
V_\lambda\big\}.$$
The maps $p_1:{\cal V}\to \D$ and $p_2:{\cal V}\to \C$ are the
projections to the first and the second coordinates.

\begin{proposition}
There exists a family of simply connected open sets $\wt{V}_\lambda$ and of
universal coverings $\pi_\lambda : \wt{V}_\lambda \to V_\lambda$ such
that $\wt{V_0} = \D$, the set
$$\wt{\cal V} = \big\{ (\lambda,z) \in \D \times\C\, \big| z \in
\wt{V}_\lambda \big\}$$ is open,
and $\Pi : (\lambda,z) \in \wt{\cal V} \mapsto \pi_\lambda(z)$ is
analytic.\\
For all $\lambda \in \D$,
$$\wt{V}_\lambda \subset B(0,\rho) \text{ with }
\log \rho = \frac{2\, \log 4}{\ds 1+ |\lambda|^{-1}}.$$
\end{proposition}
\proof
We want to construct universal coverings 
$\pi_\lambda:\widetilde V_\lambda\to V_\lambda$ such that $\pi_\lambda$ 
depend holomorphically on $\lambda$. For this purpose, we use Bers's 
embedding.

By hypothesis, the set $V_0$ is hyperbolic, i.e. its analytic
universal coverings are isomorphic to $\D$.
Let $\pi_0:\D\to V_0$ be a universal covering mapping $0$ to $0$.
Let 
$h_\lambda:V_0\to V_\lambda$ be the quasiconformal homeomorphism 
provided by the holomorphic motion. Let $\mu_\lambda$ be the Beltrami 
form on $V_0$ defined by 
$\mu_\lambda = \overline \partial h_\lambda/\partial h_\lambda$. 
Finally, let $\widetilde\mu_\lambda$ be the Beltrami form defined on $\C$ by 
$\widetilde\mu_\lambda = \pi_0^* \mu_\lambda$ on $\D$ and 
$\widetilde\mu_\lambda=0$ on $\C\setminus\D$. 

There exist quasiconformal 
homeomorphisms $\widetilde h_\lambda:\C\to\C$ such that 
$\widetilde \mu_\lambda = \overline \partial \widetilde h_\lambda/\partial 
\widetilde h_\lambda$. Those homeomorphisms are univalent outside $\D$. 
We can normalize them by the conditions 
$\widetilde h_\lambda(0)=0$ and $\widetilde h_\lambda(z)=z+{\cal O}(1)$ as 
$z\to \infty$. Then $\wt{h}_\lambda$ is uniquely defined.

Now, set $\widetilde V_\lambda = \widetilde h_\lambda(\D)$ and define
$\pi_\lambda:\widetilde V_\lambda\to V_\lambda$ by
$\pi_\lambda = h_\lambda\circ \pi_0\circ \widetilde h_\lambda^{-1}$.
Then, $\pi_\lambda$ are universal coverings. The computation to prove that
$(\lambda,z) \mapsto \pi_\lambda(z)$ is analytic is becoming well
known, but since we know no reference for this, we include the proof
here~: indeed, for every fixed $\lambda$,
the null Beltrami differential is mapped by $h_\lambda^{-1}$ to
$\wt{\mu}_\lambda$, which is mapped by $\pi_0$ to $\mu_\lambda$, and
then by $h_\lambda$ to $0$. Thus each $\pi_\lambda$ is a holomorphic function.
Then, 
$\wt \mu_\lambda$ depends holomorphically on $\lambda$, and thus,
$\partial \wt h_\lambda/\partial \overline\lambda=0$.
Thus, if take the derivative of the expression $\pi_\lambda\circ \wt h_\lambda
=h_\lambda\circ \pi_0$ with respect to $\overline \lambda$ 
in the sense of distributions, we get
$$\frac{\partial \pi_\lambda}{\partial
\overline\lambda}\Big|_{h_\lambda(z)}+\frac{\partial
\pi_\lambda}{\partial z}\Big|_{h_\lambda(z)}\cdot\frac{\partial
\wt h_\lambda}{\partial \overline \lambda}\Big|_{z}+
\frac{\partial \pi_\lambda}{\partial \overline
z}\Big|_{h_\lambda(z)}\cdot\overline{\frac{\partial \wt
h_\lambda}{\partial \lambda}}\Big|_{z} = \frac{\partial
h_\lambda}{\partial
\overline\lambda}\Big|_{\pi_0(z)}.$$
Since $\partial \wt h_\lambda/\partial \overline\lambda=0$,
$\partial \pi_\lambda/\partial \overline z=0$ and
$\partial h_\lambda/\partial \overline\lambda = 0$, the previous
expression
simplifies to
$$\frac{\partial \pi_\lambda}{\partial
\overline\lambda}\Big|_{h_\lambda(z)}=0.$$
So, by Weil's lemma, $\pi_\lambda$ depends analytically on $\lambda$.

We can now estimate the conformal radius of the sets $\widetilde V_\lambda$. 
For this purpose, note that by the area theorem, since $\tilde h_\lambda$ is 
univalent outside $\D$ and normalized to be tangent to the identity at 
$\infty$, the set $\widetilde V_\lambda$ is contained in the disk $B(0,4)$. 
The boundary moves holomorphically in $B(0,4)\setminus \{0\}$. For 
$\lambda=0$ the boundary is the unit circle. It follows from Schwarz's 
lemma that the hyperbolic distance in  $B(0,4)\setminus \{0\}$ between 
the boundary of $\widetilde V_\lambda$ and $S^1$ is less than or equal to the 
hyperbolic distance in $\D$ between $\lambda$ and $0$.
This and an elementary computation yield
$$\widetilde V_\lambda\subset B(0,\rho_2)\quad{\rm with}\quad
\log \rho_2 = \frac{\log 16}{1+|\lambda|^{-1}}.$$
\qedprop

Let us recall that $\a \in D_i=B(p_{n_i}/q_{n_i},1/q_{n_i}^2)$.
Let $r=1/q_{n_i}^2$. The real number $\a_0$ belongs to $B_i =
B(p_{n_i}/q_{n_i},1/q_{n_i}^3)$ thus
$\ds \frac{\a_0-p_{n_i}/q_{n_i}}{r}$ belongs to $\D$.
Let us apply the previous proposition to our problem with
$$\lambda=\lambda(\a)=\zeta\Big(\frac{\a-p_{n_i}/q_{n_i}}{r}\Big)$$
where $\zeta$ is any automorphism of $\D$ that
maps $\ds \frac{\a_0-p_{n_i}/q_{n_i}}{r}$ to $0$, and with
$V_\lambda=V_{i-1}(\a)$. Let
$$\phi(\delta) = \big(\lambda \circ \a(\delta),
  \chi(\delta)\big).$$
where $\a(\delta) = \frac{p_{n_i}}{q_{n_i}}+\delta^q_{n_i}$ and $\chi
= \chi_{p_{n_i}/q_{n_i}}$.
In section~\ref{sect5}, we proved that $|\a_0 - p_{n_i}/q_{n_i}|
\leq 1/2q_{n_i}^3$. Thus
$$\lambda(B_i) \subset B\Big(0, \frac{3}{2q_{n_i}} \Big).$$
This proves proposition~\ref{movingrange}.

Now, set 
$$\widetilde {\cal V}=\{(\lambda,z)~|~\lambda\in \D
~{\rm and}~z\in \widetilde V_\lambda\}.$$
We keep the notation $p_1$ and $p_2$ for the projections on the first 
and the second coordinates. 
We can lift the map $\phi:U_i\to  {\cal V}$ to a map $\tilde \phi:U_i\to 
\widetilde {\cal V}$ such that for all $\delta\in U_i$, we have 
\begin{itemize}
\item $p_1\circ \phi(\delta) = p_1\circ \tilde \phi(\delta)=\lambda$,
\item $p_2\circ \phi(\delta) = \pi_\lambda\circ p_2\circ 
\tilde\phi(\delta)$ and
\item $p_2\circ \tilde\phi(0)=0$.
\end{itemize}
We then define $\wh{\phi} = p_2 \circ \tilde{\phi}$.

\section{Acknowledgments.}

We wish to express our gratitude to A. Douady,
C. Henriksen and C.T. McMullen, for helpful discussions.

\appendix

\section{Arithmetic conventions}\label{app_cf}

By convention,
$$[a_0,a_1,a_2,\ldots] = a_0 + \frac{1}{\ds a_1+ \frac{1}{a_2+\ddots}}.$$
The $n$-th approximant of an irrational number $\a = [a_0,a_1,\ldots]
\in \R \setminus \Q$ is the
number
$$p_n/q_n = [a_0, \ldots, a_n] = a_0 + \frac{1}{\ds \ddots+\frac{1}{a_n}},$$
where $q_n$ is a positive integer, and the fraction $p_n/q_n$ is in
it's lowest terms.

We always have $q_0 = 1$, and if we set $q_{-1}=0$ and $q_{-2} =1$,
then the following recurrence relation holds for all $n \in \N$:
$$q_n = a_n q_{n-1} + q_{n-2}.$$
Thus $q_1 = a_1$, $q_2 = a_2 a_1 +1$, \ldots, and $q_n$ never depends
on $a_0$.

For $\a \in \R\setminus\Q$, the Bruno sum is defined by
$$B(\a) = \sum_{n=0}^{+\infty} \frac{\log q_{n+1}}{q_n}\in\ ]0,+\infty].$$

\begin{lemma}
For $\a \in \R\setminus\Q$, $B(\a+1) = B(\a)$ and $ B(1-\a) = B(\a)$.
\end{lemma}
\proof
The first comes from $\a+1 = [a_0+1,a_1,a_2,\ldots]$.
For the second, we may assume that $\a \in ]1/2,1[$. This is equivalent
to $a_0=0$ and $a_1 = 1$. Thus $q_1 = 1$ and $q_2 = a_2 +1$.
It is easy to check that
$1-\a = [0,a_2+1,a_3,a_4,\ldots]$.
Thus, if we note $p'_n/q'_n$ the approximants of $1-\a$, then
$q'_0 = 1 =q_1$, $q'_1=a_2+1 =q_2$, and one then checks by induction
that $q'_n = q_{n+1}$ for all $n\in\N$. Thus $B(\a) =
B(1-\a) + \log(q_1)/q_0$, and $\log(q_1)/q_0 = \log(1)/1 = 0$.
\qedprop

At some point, we defined the Fibonacci numbers $F_n$, by
$F_0=1$, $F_1=2$, and $F_{n+2} = F_{n+1} + F_n$ for all $n\in\N$.
The reader should note that the indexing may be different than what
is usually found in the litterature.
It is desinged for the situation when $\a \in\, ]0,1/2[$.
Then, for all $n \in \N$, $q_n \geq F_n$, as can be proved by induction.

In this article, we make use a few times of the following fact, that
we state here (the proof is left as an exercise to the reader)
\begin{lemma}\label{dulllemma}
 For all $\lambda > 81/64$, the sequence $$\log(\lambda n^2)/n,$$
 defined for $n \geq 2$, is decreasing. As a corollary,
 if $\a \in ]0,1/2[$ is irrationnal, then for all $n \geq 1$,
 $$ \frac{\log(\lambda q_n^2)}{q_n} \leq
 \frac{\log(\lambda F_n^2)}{F_n}.$$
\end{lemma}

We also include here, for reference, the following computations:
$$
 \sum_{n=1}^{+\infty} \frac{\log F_n}{F_n}  =  1.96\ldots \qquad
 \sum_{n=1}^{+\infty} \frac{1}{F_n}  =  1.35\ldots
$$
where the rounding is to the lower.

\newcounter{nom}{\setcounter{nom}{1}}

\end{document}